\documentclass[onecolumn,showpacs,showkeys,preprintnumbers]{revtex4}

\usepackage{amsmath}
\usepackage{amssymb}
\usepackage{graphicx}
\usepackage{dcolumn}
\usepackage{bm}
\usepackage{mathrsfs}
\usepackage{enumerate}
\usepackage{hyperref}

\newtheorem{theorem}{Theorem}[section]
\newtheorem{lemma}[theorem]{Lemma}

\DeclareMathOperator{\arccot}{arccot}%
\DeclareMathOperator{\sh}{sh}%
\DeclareMathOperator{\ch}{ch}%
\DeclareMathOperator{\cth}{cth}%
\DeclareMathOperator{\Arcth}{Ar~cth}%
\DeclareMathOperator{\Arth}{Ar~th}%
\DeclareMathOperator{\Arsh}{Ar~sh}%
\DeclareMathOperator{\Arch}{Ar~ch}%
\DeclareMathOperator{\arcsec}{arcsec}%
\DeclareMathAlphabet{\mathpzc}{OT1}{pzc}{m}{it}

\begin{document}
\title{Energy Optimal Control for Quantum System Evolving on
$SU(1,1)$}

\author{Jian-Wu Wu$^{1}$}
\email{wujw03@mails.tsinghua.edu.cn}
\author{Chun-Wen Li$^{1}$}
\author{Jing Zhang$^{1}$}
\author {Tzyh-Jong Tarn$^{2}$}

\affiliation{$^{1}$Department of Automation, Tsinghua University,
Beijing, 100084, P. R. China\\
$^{2}$Department of Systems Science and Mathematics, Washington
University, St. Louis, MO 63130, USA}

\date{\today}

\begin{abstract}
This paper discusses the energy optimal control problem for the
class of quantum systems that possess dynamical symmetry of
$SU(1,1)$, which are widely studied in various physical problems
in the quantum theory. Based on the maximum principle on Lie
group, the complete set of optimal controls are analytically
obtained, including both normal and abnormal extremals. The
results indicate that the normal extremal controls can be
expressed by the Weierstrass elliptic function, while the abnormal
extremal controls can only be constant functions of time $t$.
\end{abstract}



\maketitle

\section{Introduction}\label{sec1}

During the past two decades, optimization techniques have been
extensively applied to design external control fields to
manipulate the evolution of quantum mechanical systems~
\cite{Peirce1,Shi1,Zhu1,Schirmer1,Zhao1,khaneja1,wrb1,Alessandro2,Carlini1,Boscain2,Boscain3,Alessandro1,Boscain1}.
A typical application is to force the states to approach the
priori prescribed targets as closely as possible, including both
bounded~\cite{Peirce1,Shi1,Zhu1,Schirmer1} and
unbounded~\cite{Zhao1} situations. Also, optimization theory can
be applied to improve the efficiency of desired quantum state
transitions, e.g., the evolution time
\cite{khaneja1,wrb1,Alessandro2,Carlini1,Boscain2,Boscain3} and
the energy consumed by the controls ~\cite{Boscain1,Alessandro1}
are most interesting. For some low dimensional quantum systems
that evolve on compact Lie groups, one can find analytical
solutions for such optimization problems with
bounded~\cite{Alessandro1,wrb1} or unbounded~\cite{khaneja1}
controls. However, in the higher dimensional situations, numerical
algorithms have to be applied.

In this paper, we explore the optimal steering problem for the
class of quantum systems whose evolution operators can be
described by the $SU(1,1)$ matrices. The underlying system is
modelled by an evolution equation of the form~\cite{wjw1}:%
\begin{equation}\label{eq1.1}
\frac{d}{dt}X(t)=\left[A+u(t)B\right]X(t),~X(0)=I_2,
\end{equation}

\noindent where $X(t)$ is a two dimensional special pseudo-unitary
matrix; $u(t)$ is real function of time $t$, which is the control
input of the system; $A$ and $B$ are arbitrary matrices that can
be expressed as the linear combination of $K_x$, $K_y$ and $K_z$,
which are the generators of the Lie algebra $su(1,1)$ and can be
identified as follows:%
\begin{equation}\label{eq1.1a}
K_x=\frac{1}{2}\left(
\begin{array}{cc}
  0 & -i \\
  i & 0\\
\end{array}\right),~~
K_y=\frac{1}{2}\left(
\begin{array}{cc}
  0 & -1 \\
  -1 & 0 \\
\end{array}\right),~~
K_z=\frac{1}{2}\left(
\begin{array}{cc}
  -i & 0 \\
  0 & i \\
\end{array}\right).
\end{equation}

\noindent In \cite{Jurdjevic4}, Jurdjevic has taken the initial
steps in the problems of optimal control for the special case when
$A=K_x$ and $B=K_z$~(the corresponding properties of the optimal
controls also can be found for the systems evolving on the
homomorphic groups $SO(2,1)$ and $SL(2,\mathbb{R})$ in
\cite{Jurdjevic4} and \cite{Jurdjevic3}). Unlike the case of
$SU(2)$~\cite{Alessandro3}, however, as will be seen in
Section~\ref{sec2}, system (\ref{eq1.1}) usually can't be
transformed into the special form such that $A=K_x$ and $B=K_z$.
Thus, it is natural for us to consider the general case in detail.

With specific realizations and representations of the Lie algebra
$su(1,1)$ introduced, system (\ref{eq1.1}) can be used to describe
various of quantum dynamical processes, e.g., the superfluid
system under Bose realization~\cite{Solomon1}, the harmonic
oscillator under $xp$-realization~\cite{Penna1}, the $SU(1,1)$
coherent states under irreducible unitary representation with
respect to positive discrete series~\cite{Gerry1}.

The optimal control problem to be considered in this paper is
formulated as follows. Given an arbitrary target evolution matrix
$X_f$ in $SU(1,1)$, we wish to find a control function $u(t)$ that
can steer the evolution matrix associated with system
(\ref{eq1.1}) from its initial state $I_2$ to some desired final
state $X_f$, and meanwhile, minimize the quadratic cost function%
\begin{equation}\label{eq1.2}
J(u)=\int_0^Tu(t)^2dt,
\end{equation}

\noindent where $T$ is the final time. The quadratic cost index
given in (\ref{eq1.2}) measures the energy consumed during the
steering process between the initial $I_2$ and the terminal $X_f$.
Since the remarkable difference between the quantum systems and
the classical systems is that the evolutions of the former may be
disturbed by the decoherence phenomenon. This is practical because
that the increasing of the intensity of the electromagnetic
fields, which are used to control the evolution of the coupled
quantum system, tend to induce relaxation and decoherence
phenomena.

Based on the maximum principle for systems evolving on Lie
group~\cite{Baillieul1}, explicit forms of the control functions
with respect to both normal and abnormal extremals will be derived
analytically. The problem considered here can be viewed as the
noncompact prolongation of the $SU(2)$ case presented
in~\cite{Alessandro1}. However, according to quantum theory and
group representation theory~\cite{Vilenkin1}, the noncompact
$SU(1,1)$ Lie group has only infinite dimensional unitary
representations, and hance the associated evolution operator (or
propagator) corresponding to $X(t)$ is always infinite
dimensional, which describes the transition between two quantum
states defined in an infinite dimensional Hilbert space. The
dynamics of the $SU(1,1)$ coherent states is a typical
example~\cite{Gerry1}. This makes the derivation quite different
from and far more complicated than that in the case of $SU(2)$.

The balance of this paper is organized as follows. In
Section~\ref{sec2}, some useful results including controllability
properties and the maximum principle with respect to the quantum
system evolving on Lie group $SU(1,1)$ are introduced. In
Section~\ref{sec3}, we discuss the optimal steering problem with
respect to the abnormal extremals. Properties of the abnormal
optimal control function are characterized. In Section~\ref{sec4},
the control functions corresponding to the normal extremals are
derived analytically for all possible cases. In
Section~\ref{sec5}, two examples are provided for illustration.
Finally, conclusions are drawn in Section~\ref{sec6}.

\section{Preliminaries on the Quantum Control System on the Lie Group $SU(1,1)$}\label{sec2}

With the three generators $K_x$, $K_y$ and $K_z$ of the Lie
algebra $su(1,1)$ in (\ref{eq1.1a}), any given evolution matrix
$X$ associated with system
(\ref{eq1.1}) can be written as%
\begin{equation}\label{eq2.2}
X=e^{{\alpha}K_z}e^{{\beta}K_y}e^{{\gamma}K_z},
\end{equation}

\noindent where $-2\pi<\alpha,\gamma\leq2\pi$,
$0\leq\beta<\infty$. The commutation relations between $K_x$,
$K_y$ and $K_z$ are%
\begin{equation}\label{eq2.3}
[K_x,K_y]=-K_z,~~~[K_y,K_z]=K_x,~~~[K_z,K_x]=K_y.
\end{equation}

\noindent With the inner product $\left<\cdot,\cdot\right>$ defined by%
\begin{equation}\label{eq2.4}
\left<M,N\right>=2\mbox{Tr}(MN^\dag),
\end{equation}

\noindent where $N^\dag$ is the Hermitian conjugation of $N$, it
can be verified that $K_x$, $K_y$ and $K_z$ form an orthonormal
basis of the Lie algebra $su(1,1)$. Accordingly, the drift term
$A$ and the control term $B$ in system (\ref{eq1.1}) can be
expressed by linear combinations of the three generators $K_x$,
$K_y$ and $K_z$ as follows:%
\begin{eqnarray}\label{eq2.5}
 A &=&\left<A,K_x\right>K_x+\left<A,K_y\right>K_y+\left<A,K_z\right>K_z,\\\label{eq2.6}
 B &=&\left<B,K_x\right>K_x+\left<B,K_y\right>K_y+\left<B,K_z\right>K_z.
\end{eqnarray}

In this paper, we say $M$ is pseudo-orthogonal to $N$ when
$\left<M,N^\dag\right>=0$. And $M$ will be called elliptic
(hyperbolic, parabolic) if $\left<M,M^\dag\right>$ is negative
(positive, zero). Accordingly, the matrices in the Lie algebra
$su(1,1)$ are separated into three different types.  Since for any
matrix $P{\in}SU(1,1)$,
we have%
\noindent\begin{equation}\label{eq2.5a}
\left<PMP^{-1},(PMP^{-1})^\dag\right>=2\mbox{Tr}(PMP^{-1}PMP^{-1})=2\mbox{Tr}(MM)=\left<M,M^\dag\right>,
\end{equation}

\noindent i.e., none of the changes of coordinates by the
$SU(1,1)$ transformations will alter the type of an $su(1,1)$
matrix. Thus, as mentioned in Section~\ref{sec1}, system
(\ref{eq1.1}) usually can't be transformed into the special case
such that $A=K_x$~(hyperbolic) and $B=K_z$~(elliptic).

We assume system (\ref{eq1.1}) is controllable in this paper, this
ensures that the optimal control problem stated in the
introduction section is always solvable. System (\ref{eq1.1}) is
said to be controllable on $SU(1,1)$ if for any given target
evolution matrix $X_f{\in}SU(1,1)$ there always exists at least
one control $u(t)$ such that $X(u;T)=X_f$ for some time $T$. The
problems of the controllability for systems evolving on Lie groups
have received a great deal of attention in the past decades (see
e.g., \cite{Jurdjevic1,Jurdjevic2}). As for the systems evolving
on the noncompact Lie group $SU(1,1)$, a sufficient and necessary
condition is provided in~\cite{wjw2} and can be summarized as
follows.%

\begin{theorem}\label{theorem2.1}%
(\cite{wjw2})
\begin{enumerate}
    \item If $A$ and $B$ are linearly dependent, then system (\ref{eq1.1}) is uncontrollable on
    $SU(1,1)$;
    \item If $A$ and $B$ are linearly independent, then system (\ref{eq1.1}) is controllable if and
    only if the set $\{u\in\mathbb{R}|\left<A\right.+uB,A^{\dag}+\left.uB^{\dag}\right><0\}$
    is nonempty.
\end{enumerate}
\end{theorem}%

To obtain the optimal control functions, we will make use of the
well-known maximum principle on Lie
groups~\cite{Baillieul1,Jurdjevic2}, which, under the assumption
that system (\ref{eq1.1}) is controllable, can be summarized for
the $SU(1,1)$ case as follows~\cite{Baillieul1}.

\noindent\begin{theorem}\label{theorem2.2}%
Assume that system (\ref{eq1.1}) is controllable, if $u^o(t)$ is
an optimal control that minimizes the quadratic index given in
(\ref{eq1.2}) and $X^o(t)$ is the corresponding optimal trajectory
of (\ref{eq1.1}). Then, there exists a constant matrix
$S{\in}su(1,1)$ and a nonnegative real number $\lambda_o$, not
both zero, such that for almost everywhere (a.e.) the Hamiltonian function%
\begin{equation}\label{eq2.19}
H(S;\lambda_o;u;X_o(t))=\left<S,{X^o}(t)^{-1}[A+u(t)B]{X^o}(t)\right>+\frac{1}{2}\lambda_ou(t)^2
\end{equation}

\noindent is minimized with respect to $u$ by $u^o(t)$.
\end{theorem}

The above theory immediately provides a necessary condition for
the optimality. The problems with respect to $\lambda_o\neq0$ are
called normal, otherwise are called abnormal.

\section{Optimal Control Function for the Abnormal Extremal}\label{sec3}

In this section, we will investigate the abnormal case. For this
purpose, the following basic properties on $su(1,1)$ are useful.

\begin{lemma}\label{lemma3.1}%
For arbitrary pair of matrices $M$ and $N$ in $su(1,1)$, the
following relations hold:%
\begin{eqnarray}\label{eq3.1}
  [[M,N],M]=\left<M,N^{\dag}\right>M-\left<M,M^{\dag}\right>N,\\\label{eq3.2}
  [[M,N],N]=\left<N,N^{\dag}\right>M-\left<M,N^{\dag}\right>N.~
\end{eqnarray}
\end{lemma}

{\it{proof}} Since (\ref{eq3.1}) is equivalent to (\ref{eq3.2})
because of symmetry, we only need to proof (\ref{eq3.1}). If $M$
and $N$ are linearly dependent, then
there exists a real constant $k\neq0$ such that $M=kN$. Therefore,%
\begin{equation}\label{eq3.3}%
[[M,N],M]=[[kN,N],kN]=0,
\end{equation}%

\noindent and%
\begin{equation}\label{eq3.4}
\left<M,N^{\dag}\right>M-\left<M,M^{\dag}\right>N=\left<kN,N^{\dag}\right>kN-\left<kN,kN^{\dag}\right>N=0,
\end{equation}%

\noindent i.e., the relation in (\ref{eq3.1}) holds. Thus, we only
need to consider the case when $M$ and $N$ are linearly
independent, i.e., $[M,N]\neq0$. Since%
\begin{equation}\label{eq3.6}
\left<[M,N]^{\dag},M\right>=\left<[M,N],M^{\dag}\right>=\left<MN,M^{\dag}\right>-\left<NM,M^{\dag}\right>=0,
\end{equation}

\noindent and%
\begin{equation}\label{eq3.7}
\left<[M,N]^{\dag},N\right>=\left<[M,N],N^{\dag}\right>=\left<MN,N^{\dag}\right>-\left<NM,N^{\dag}\right>=0,
\end{equation}
$M$, $N$ and $[M,N]^{\dag}$ are still linearly independent. Thus
we can express $[[M,N],~M]$ as%
\begin{equation}\label{eq3.8}
[[M,N],M]={\mu_1}M+{\mu_2}N+{\mu_3}[M,N]^{\dag},
\end{equation}

\noindent for some constants $\mu_1$, $\mu_2$ and $\mu_3$. Taking
the inner product with $[M,N]^{\dag}$ in (\ref{eq3.8}), we obtain:%
\begin{equation}\label{eq3.10}
{\mu_3}\left<[M,N]^{\dag},[M,N]^{\dag}\right>=0.
\end{equation}

\noindent Eq.~(\ref{eq3.10}) implys that ${\mu_3}=0$, thus we can
rewrite (\ref{eq3.8}) as:%
\begin{equation}\label{eq3.12}
[[M,N],M]={\mu_1}M+{\mu_2}N.
\end{equation}

\noindent Taking the inner product of (\ref{eq3.12}) with
$M^{\dag}$, we get%
\begin{equation}\label{eq3.14}
{\mu_1}\left<M,M^{\dag}\right>+{\mu_2}\left<N,M^{\dag}\right>=0.~
\end{equation}

\noindent Since%
\begin{equation}\label{eq3.16}
\left<[[M,N],M],N^{\dag}\right>=2\mbox{Tr}([[M,N],M]N)=2\mbox{Tr}([M,N][M,N])=\left<[M,N],[M,N]^{\dag}\right>,
\end{equation}

\noindent take the inner product of (\ref{eq3.12}) with
$N^{\dag}$, we have%
\begin{equation}\label{eq3.18}
\left<[M,N],[M,N]^{\dag}\right>={\mu_1}\left<M,N^{\dag}\right>+{\mu_2}\left<N,N^{\dag}\right>.
\end{equation}

\noindent Make use of the equality
$\left<[M,N],[M,N]^{\dag}\right>={\left<M,N^{\dag}\right>^2-\left<M,M^{\dag}\right>\left<N,N^{\dag}\right>}$~(see
the proof in \cite{wjw2}), from (\ref{eq3.14}) and (\ref{eq3.18})
we have%
\begin{equation}\nonumber
  {\mu_1}=\left<M,N^{\dag}\right>,~~~~%
  {\mu_2}=-\left<M,M^{\dag}\right>.%
\end{equation}

\noindent\begin{lemma}\label{lemma3.2}%
If $M$ and $N$ are linearly independent and the set
$\{u\in\mathbb{R}|\left<M\right.+uN,M^{\dag}+\left.uN^{\dag}\right><0\}$
is nonempty, then $M$, $N$ and $[M,N]$ form a basis in $su(1,1)$.
\end{lemma}

\noindent{\it{proof}} Assume that $[M,N]$ can be linearly
expressed by $M$ and $N$, i.e. there exist two real number
$\lambda_1$ and $\lambda_2$ satisfy%
\begin{equation}\label{eq3.20}
[M,N]=\lambda_1{M}+\lambda_2{N}.
\end{equation}

\noindent On the one hand, making communications with $M$ and $N$
respectively, one can obtain%
\begin{eqnarray}\label{eq3.21}
[[M,N],M]=-\lambda_2[M,N]=-\lambda_1\lambda_2{M}-{\lambda_2}^2{N},\\\label{eq3.22}
[[M,N],N]=\lambda_1[M,N]={\lambda_1}^2{M}+\lambda_1\lambda_2{N}.~~
\end{eqnarray}

\noindent On the other hand, notice that (\ref{eq3.1}) and
(\ref{eq3.2}) would imply%
\begin{eqnarray}\label{eq3.23}
  [[M,N],M]=\left<M,N^{\dag}\right>M-\left<M,M^{\dag}\right>N,\\\label{eq3.24}
  [[M,N],N]=\left<N,N^{\dag}\right>M-\left<M,N^{\dag}\right>N.~
\end{eqnarray}

\noindent Compare the coefficients of (\ref{eq3.21}) and
(\ref{eq3.22}) and those of (\ref{eq3.23}) and (\ref{eq3.24})
respectively, we obtain%
\begin{equation}\label{eq3.25}
{\lambda_1}^2=\left<N,N^{\dag}\right>,~~
\lambda_1\lambda_2=-\left<M,N^{\dag}\right>,~~
{\lambda_2}^2=\left<M,M^{\dag}\right>.
\end{equation}

\noindent A straightforward computation, using (\ref{eq3.25}),
shows that, for every $u\in\mathbb{R}$,%
\begin{eqnarray}\nonumber
&&\left<M+uN,M^{\dag}+uN^{\dag}\right>\\\nonumber
&=&\left<M,M^{\dag}\right>+2u\left<M,N^{\dag}\right>+u^2\left<N,N^{\dag}\right>\\\nonumber
&=&{\lambda_2}^2-2\lambda_1\lambda_2u+{\lambda_1}^2u^2\\\label{eq3.26}
&=&({\lambda_2}-{\lambda_1}u)^2\geq0.
\end{eqnarray}

\noindent This contradicts with the condition that the set
$\{u\in\mathbb{R}|\left<M+uN,M^{\dag}+uN^{\dag}\right><0\}$ is
nonempty. So $M$, $N$ and $[M,N]$ are linearly independent.

With the properties on the Lie algebra $su(1,1)$ obtained above,
we can now draw the conclusion for the abnormal case as follows.

\begin{theorem}\label{theorem3.3}%
Assume that the system (\ref{eq1.1}) is controllable. If
$\left<B,B^{\dag}\right>\neq0$, then the control function
$u(t)=-\left<A,B^{\dag}\right>/\left<B,B^{\dag}\right>$, a.e., is
the only abnormal extremal for the optimization problem under
consideration. Otherwise, there is no abnormal extremal.
\end{theorem}

{\it{proof}}
If $\lambda_o=0$, the Hamiltonian function may be rewritten as%
\begin{equation}\label{eq3.27}
H(S;\lambda_o;u;X_o(t))=\left<S,{X^o}(t)^{-1}[A+u(t)B]{X^o}(t)\right>.
\end{equation}

\noindent Since the minimization condition indicates that
$H(S;\lambda_o;u;X_o(t))$ is a.e. minimized with respect to $u(t)$
by $u_o(t)$. To obtain the optimal control function $u_o(t)$, we
differentiate $H(S;\lambda_o;u;X_o(t))$ with respect ro $u$ and
set the result equal to zero, which yields%
\begin{equation}\label{eq3.28}
\left<S,{X^o}(t)^{-1}B{X^o}(t)\right>=0,~~\text{a.e}.
\end{equation}

\noindent By continuity, (\ref{eq3.28}) implies that%
\begin{equation}\label{eq3.28A}
\left<S,{X^o}(t)^{-1}B{X^o}(t)\right>\equiv0.
\end{equation}

\noindent Differentiating (\ref{eq3.28A}) with respect ro $t$,
using (\ref{eq1.1}), we can obtain%
\begin{equation}\label{eq3.29}
\left<S,{X^o}(t)^{-1}[A,B]{X^o}(t)\right>=0.
\end{equation}%
\noindent Differentiating (\ref{eq3.29}) with respect ro $t$ again, we have%
\begin{equation}\label{eq3.30}
\left<S,{X^o}(t)^{-1}[[A,B],A]{X^o}(t)\right>+u(t)\left<S,{X^o}(t)^{-1}[[A,B],B]{X^o}(t)\right>=0.
\end{equation}

\noindent Utilizing  (\ref{eq3.1}), (\ref{eq3.2}) and
(\ref{eq3.28}), we can recast (\ref{eq3.30}) to%
\begin{equation}\label{eq3.31}
\left(\left<A,B^{\dag}\right>+u(t)\left<B,B^{\dag}\right>\right)\left<S,{X^o}(t)^{-1}A{X^o}(t)\right>=0,~~\text{a.e}.
\end{equation}

\noindent Notice that, according to Lemma\ref{lemma3.2}, $A$, $B$
and $[A,B]$ are linearly independent. If the equality
$\left<S,\right.{X^o}(t)^{-1}A\left.{X^o}(t)\right>=0$ holds, then
combining (\ref{eq3.28}) and (\ref{eq3.29}) we can draw conclusion
that $S$ must be zero, which is a contradiction. Therefore
Eq.(\ref{eq3.31})
implies that%
\begin{equation}\label{eq3.32}
\left<A,B^{\dag}\right>+u(t)\left<B,B^{\dag}\right>=0,~~\text{a.e}.
\end{equation}

\noindent If $\left<B,B^{\dag}\right>=0$, from Eq.(\ref{eq3.32})
we have $\left<A,B^{\dag}\right>=0$, which contradicts the
assumption that system (\ref{eq1.1}) is controllable~(see
\cite{wjw2}). Thus, there is no abnormal extremal when
$\left<B,B^{\dag}\right>=0$. If $\left<B,B^{\dag}\right>\neq0$,
then from (\ref{eq3.32}) we have
$u(t)=-{\left<A,~B^{\dag}\right>}/{\left<B,~B^{\dag}\right>}$,
a.e.

{\it{Remark:}} The target evolution matrices that can be achieved
directly by the abnormal control
$u=-{\left<A,~B^{\dag}\right>}/{\left<B,~B^{\dag}\right>}$ are in
the one dimensional Lie group corresponding to the Lie subalgebra
of $su(1,1)$ generated by
$A-{\left<A,~B^{\dag}\right>}/{\left<B,~B^{\dag}\right>}B$, which
can never fill up the whole Lie group of $SU(1,1)$. Consequently,
abnormal extremal exists only when the target evolution matrix
$X_f$ is of the form
${\exp}\left[c\left(A-{\left<A,~B^{\dag}\right>}/{\left<B,~B^{\dag}\right>}B\right)\right]$
for some real constant $c$. In order to solve the optimal steering
problem subject to the terminal condition $X(T)=X_f$, where $X_f$
is an arbitrary matrix taken from $SU(1,1)$, the candidates can
only be normal extremals.

\section{Optimal Control Function for the Normal Extremal}\label{sec4}

In this section we explore the normal extremal control functions.
According to Theorem~\ref{theorem2.1}, here we assume that $A$ and
$B$ are linearly independent and meanwhile the set
$\{u\in\mathbb{R}|\left<A\right.+uB,A^{\dag}\left.+uB^{\dag}\right><0\}$
is nonempty throughout this section, to guarantee the
controllability of the system.

One can obtain the normal extremal, according to
Theorem~\ref{theorem2.2}, by minimizing the Hamiltonian function
$H(S;\lambda_o;u;X_o(t))$ as a quadratic function of $u$. After
normalizing $\lambda_o=1$, by continuity, the necessary condition
for candidate optimal controls can be expressed as%
\begin{equation}\label{eq4.1}
    u(t)=-\left<S,X_o^{-1}(t)BX_o(t)\right>, \text{a.e.}
\end{equation}

\noindent where the matrix $S$ is an element in the Lie algebra
$su(1,1)$. We introduce the following two auxiliary variables in
the succeeding discussion:%
\begin{equation}\label{eq4.2}
\begin{array}{l}
  u_A(t)=-\left<S,{X}_o(t)^{-1}AX_o(t)\right>,\\
  u_C(t)=-\left<S,{X}_o(t)^{-1}[A,B]X_o(t)\right>.\\
\end{array}
\end{equation}

\noindent In order to follow standard notations, we rewrite the
normal extremal $u(t)$ in (\ref{eq4.1}) as $u_B(t)$. Differentiate
$u_A$, $u_B$ and $u_C$ with respect to the time $t$, respectively,
and make use of (\ref{eq1.1}),
(\ref{eq3.1}) and (\ref{eq3.2}), we can obtain (a.e.)%
\begin{eqnarray}\label{eq4.4}
   &&\dot{u}_A=u_Bu_C,\\\label{eq4.5}
   &&\dot{u}_B=-u_C,\\\label{eq4.6}
   &&\dot{u}_C={\alpha}u_A-{\beta}u_B+{\gamma}u_Au_B-{\alpha}u_B^2,
\end{eqnarray}

\noindent where%
\begin{equation}\label{eq4.4-6}
    \alpha=\left<A,B^{\dag}\right>,~~\beta=\left<A,A^{\dag}\right>,~~\gamma=\left<B,B^{\dag}\right>.
\end{equation}

From (\ref{eq4.4})-(\ref{eq4.6}), it is easy to verify the
following conclusion for the normal extremal $u_B$ and the two
auxiliary variables $u_A$ and $u_C$.

\begin{theorem}\label{theorem4.1}
The following two quantities are conserved along the normal
extremal trajectories, i.e.,%
\begin{eqnarray}\label{eq4.7}
&&u_A+\frac{1}{2}u_B^2=c_1,\\\label{eq4.8}
&&\frac{1}{2}{\gamma}u_A^2-{\alpha}u_Au_B-\frac{1}{2}u_C^2-{\beta}u_A=c_2,
\end{eqnarray}
\noindent for some constants $c_1$ and $c_2$.
\end{theorem}

According to Theorem~\ref{theorem4.1}, the initial and the final
values of $u_A$, $u_B$ and $u_C$ should satisfy (a.e.)%
\begin{equation}\label{eq4.8a}
\begin{array}{l}
   u_A(0)+\frac{1}{2}u_B(0)^2=u_A(T)+\frac{1}{2}u_B(T)^2,\\
   \frac{1}{2}{\gamma}u_A(0)^2-{\alpha}u_A(0)u_B(0)-\frac{1}{2}u_C(0)^2-{\beta}u_A(0)\\
~~~~~~~=\frac{1}{2}{\gamma}u_A(T)^2-{\alpha}u_A(T)u_B(T)-\frac{1}{2}u_C(T)^2-{\beta}u_A(T),\\
\end{array}
\end{equation}
\noindent where%
\begin{equation}\label{eq4.3a}
    u_A(0)=-\left<S,A\right>,~~~
    u_B(0)=-\left<S,B\right>,~~~
    u_C(0)=-\left<S,[A,B]\right>,
\end{equation}

\noindent and%
\begin{equation}\label{eq4.3b}
\begin{array}{l}
  u_A(T)=-\left<S,X(T)^{-1}AX(T)\right>, \\
  u_B(T)=-\left<S,X(T)^{-1}BX(T)\right>, \\
  u_C(T)=-\left<S,X(T)^{-1}[A,B]X(T)\right>, \\
\end{array}
\end{equation}

\noindent The matrix $S$ in Eqs.(\ref{eq4.3a}) and (\ref{eq4.3b})
then can be viewed as parameter matrix, which has to be chosen to
match the final condition $X(T)=X_f$.%

\begin{center}
\begin{tabular}{@{}c@{}|@{}c@{}}
\multicolumn{2}{c}{\vspace{0.2cm}\text{Table~I. The
controllability of system~(\ref{eq1.1}) with respect to the values
of $\alpha$,
$\beta$ and $\gamma$.}}\\
 \hline\hline
The values of $\alpha$,~$\beta$,~and $\gamma$  &~~the controllability of system~(\ref{eq1.1})~~\\\hline%

\begin{tabular}{@{}c@{}|@{}c@{}}
~~~~~~~~$\alpha=0$~~~~~~~~&

\begin{tabular}{@{}c@{}}%
\begin{tabular}{@{}c@{}|@{}c@{}}%
~~~~~~~~~~$\gamma<0$~~~~~~~~~~&%
\begin{tabular}{@{}c@{}}%
~~~~~~~~~~$\beta=0$~~~~~~~~~~\\\hline%
~~~~~~~~~~$\beta\neq0$~~~~~~~~~~%
\end{tabular}%
\end{tabular}\\\hline%

$\gamma=0$\\\hline%

\begin{tabular}{@{}c@{}|@{}c@{}}%
~~~~~~~~~~$\gamma>0$~~~~~~~~~~&
\begin{tabular}{@{}c@{}}
~~~~~~~~~~$\beta<0$~~~~~~~~~~\\\hline%
~~~~~~~~~~$\beta\geq0$~~~~~~~~~~%
\end{tabular}%
\end{tabular}%

\end{tabular}\\\hline%
~~~~~~~~$\alpha\neq0$~~~~~~~~&%

\begin{tabular}{@{}c@{}}%
$\gamma\leq0$\\\hline%
\begin{tabular}{@{}c@{}|@{}c@{}}%
~~~~$\gamma>0$~~~~&
\begin{tabular}{@{}c@{}}%
$\beta\leq0$\\\hline
\begin{tabular}{@{}c@{}|@{}c@{}}
~~~~$\beta>0$~~~~&%
\begin{tabular}{@{}c@{}}%
~~~~$\alpha^2-\beta\gamma\leq0$~~~~\\\hline%
~~~~$\alpha^2-\beta\gamma>0$~~~~
\end{tabular}%
\end{tabular}%
\end{tabular}%
\end{tabular}%
\end{tabular}%
\end{tabular}%

&%

\begin{tabular}{@{}c@{}}%
~~~~~~~~~~~~~{uncontrollable}~~~~~~~~~~~~~~\\\hline%
~~~~~~~~~~~~~~{controllable}~~~~~~~~~~~~~~~\\\hline%
~~~~~~~~~~~~~{uncontrollable}~~~~~~~~~~~~~~\\\hline%
~~~~~~~~~~~~~~{controllable}~~~~~~~~~~~~~~~\\\hline%
~~~~~~~~~~~~~{uncontrollable}~~~~~~~~~~~~~~\\\hline%
~~~~~~~~~~~~~~{controllable}~~~~~~~~~~~~~~~\\\hline%
~~~~~~~~~~~~~~{controllable}~~~~~~~~~~~~~~~\\\hline%
~~~~~~~~~~~~~{uncontrollable}~~~~~~~~~~~~~~\\\hline%
~~~~~~~~~~~~~~{controllable}~~~~~~~~~~~~~~~
\end{tabular}%
\\\hline\hline%
\end{tabular}%
\end{center}

Making use of (\ref{eq4.7}) and (\ref{eq4.8}), from
(\ref{eq4.4})-(\ref{eq4.6}), we can obtain the following
differential equation for the candidate optimal control $u(t)$
given in (\ref{eq4.1})%
\begin{equation}\label{eq4.9}
(\dot{u})^2=\frac{\gamma}{4}u^4+\alpha{u^3}+(\beta-\gamma{c_1})u^2-2\alpha{c_1}u+\gamma{c_1}^2-2\beta{c_1}-2c_2,
\text{a.e.}
\end{equation}

\noindent where $u(0)=u_B(0)$. We will show, in the following,
that the candidate optimal control function can be analytically
solved from (\ref{eq4.9}) in terms of the Weierstrass
function.%

Since the involved system is assumed to be controllable to ensure
that the optimal steering problem has solutions, one only need to
consider the controllable situations accordingly. Table~I shows
the controllability properties of system (\ref{eq1.1}) in
different cases.~(see~\cite{wjw2} for details). There are three
different cases, which need to be taken into account, depending on
the values of $\alpha$ and $\gamma$.

\noindent {\it{1.~Case $\alpha=0$ and $\gamma\neq0$.}}

In this case, the drift term $A$ of system~(\ref{eq1.1}) is
pseudo-orthogonal to the control term $B$ while the latter is not
parabolic. Accordingly, Eq.(\ref{eq4.9}) can be simplified as (a.e.)%
\begin{equation}\label{eq4.10}
(\dot{u})^2=\frac{\gamma}{4}u^4+(\beta-\gamma{c_1})u^2+\gamma{c_1}^2-2\beta{c_1}-2c_2.
\end{equation}

\noindent In order to obtain the explicit form of the optimal
control function from (\ref{eq4.10}), by a variable replacement
$x=\frac{\gamma}{4}u^2+\frac{\beta-\gamma{c_1}}{3}$, we rewrite
this differential equation as (a.e.)%
\begin{equation}\label{eq4.11}
(\dot{x})^2=4x^3-g_2x-g_3,
\end{equation}

\noindent where
$g_2=\frac{1}{3}[(\beta-\gamma{c_1})^2+3(\beta^2+2\gamma{c_2})]$
and
$g_3=\frac{1}{27}[(\beta-\gamma{c_1})^2-9(\beta^2+2\gamma{c_2})]$.

From the classical theory of elliptic functions~(see, e.g.,
\cite{Whittaker1,Lawden1}), it is well known that the above
differential equation is satisfied by the Weierstrass function
$\mathscr{G}(\cdot;g_2,g_3)$ when the discriminant $g_2^3-27g_3^2$
is nonzero. Therefore, one can express the candidate
optimal control function as (a.e.)%
\begin{equation}\label{eq4.12}
u(t)=\pm\sqrt{\left.\frac{4}{\gamma}\right(\mathscr{G}(t+a;g_2,g_3)-\left.\frac{\beta-\gamma{c_1}}{3}\right)},
\end{equation}

\noindent when $g_2^3-27g_3^2\neq0$, where
$a=\mathscr{G}^{-1}\left(\frac{\gamma}{4}u_B(0)^2+\frac{\beta-\gamma{c_1}}{3};g_2,g_3\right)$.
The sign of the above candidate optimal control function $u(t)$
turns at the point when $u(t)$ cross the $t$ axis.

Consider the exceptional situations that the discriminant
of (\ref{eq4.11}) is zero, i.e.,%
\begin{equation}\label{eq4.13}
g_2^3-27g_3^2=(\beta^2+2\gamma{c_2})[(\beta^2+2\gamma{c_2})-(\beta-\gamma{c_1})^2]^2=0,
\end{equation}

\noindent It is easy to see that
either~(i)~$\beta^2+2\gamma{c_2}=0$
or~(ii)~$(\beta^2+2\gamma{c_2})-(\beta-\gamma{c_1})^2=0$.

For the case~(i), a further use of (\ref{eq4.8}) leads to that
$\beta^2+2\gamma{c_2}=(\gamma{u_A}-\beta)^2-\gamma{u_C}^2=0$.
Thus, if $\gamma<0$, we have ${u_A}=\frac{\beta}{\gamma}$,
$u_B=u_B(0)$ and $u_C=0$, which determines the candidate optimal
control function by $u(t)=u_B(0)$~(a.e.). If $\gamma>0$, from
(\ref{eq4.10}), we have~(a.e.)%
\begin{equation}\label{eq4.14}
(\dot{u})^2=\frac{\gamma}{4}[u^2+\frac{2}{\gamma}(\beta-\gamma{c_1})]^2,
\end{equation}

\noindent whose solutions can be expressed as follows.%
\begin{itemize}
    \item If $\beta-\gamma{c_1}<0$, then (a.e.)%
\begin{equation}\label{eq4.15}
    u(t)=\left\{%
          \begin{array}{ll}
             \pm\sqrt{\frac{2(\gamma{c_1}-\beta)}{\gamma}}\tanh\left(\sqrt{\frac{\gamma{c_1}-\beta}{2}}t+a\right), &
             \hbox{when $|u_B(0)|<\sqrt{\frac{2(\gamma{c_1}-\beta)}{\gamma}}$;} \\
             \pm\sqrt{\frac{2(\gamma{c_1}-\beta)}{\gamma}}\coth\left(\sqrt{\frac{\gamma{c_1}-\beta}{2}}t+a\right), &
             \hbox{when $|u_B(0)|>\sqrt{\frac{2(\gamma{c_1}-\beta)}{\gamma}}$,} \\
         \end{array}%
         \right.
\end{equation}

\noindent where%
\begin{equation}\label{eq4.16}
     a=\left\{%
          \begin{array}{ll}
          \Arth\left(\pm\sqrt{\frac{\gamma}{2(\gamma{c_1}-\beta)}}u_B(0)\right),&
          \hbox{when $|u_B(0)|<\sqrt{\frac{2(\gamma{c_1}-\beta)}{\gamma}}$;} \\
          \Arcth\left(\pm\sqrt{\frac{\gamma}{2(\gamma{c_1}-\beta)}}u_B(0)\right),&
          \hbox{when $|u_B(0)|>\sqrt{\frac{2(\gamma{c_1}-\beta)}{\gamma}}$.} \\
          \end{array}%
         \right.
\end{equation}

    \item If $\beta-\gamma{c_1}=0$, then (a.e.)%
\begin{equation}\label{eq4.17}
u(t)=\pm\frac{2}{\sqrt{\gamma}(t+a)},
\end{equation}

\noindent where
              $a=\pm\frac{\sqrt{\gamma}}{2}u_B(0)$.
    \item If $\beta-\gamma{c_1}>0$, then (a.e.)%
\begin{equation}\label{eq4.18}
    u(t)=\pm\sqrt{\frac{2(\gamma{c_1}-\beta)}{\gamma}}\tan\left(\sqrt{\frac{\gamma{c_1}-\beta}{2}}t+a\right),
\end{equation}

\noindent where
    $a=\arctan\left(\pm\sqrt{\frac{\gamma}{2(\gamma{c_1}-\beta)}}u_B(0)\right)$.
\end{itemize}

For the case~(ii), the equation (\ref{eq4.10}) can be simplified
as (a.e.)%
\begin{equation}\label{eq4.19}
(\dot{u})^2=\frac{\gamma}{4}[u^2+\frac{4}{\gamma}(\beta-\gamma{c_1})]u^2.
\end{equation}

\noindent If $\gamma<0$, clearly, the above differential equation
has no nontrivial solution other then $u(t)\equiv0$ when
$\beta-\gamma{c_1}<0$. When $\beta-\gamma{c_1}>0$, one can
immediately obtain the optimal control function
from (\ref{eq4.19}) as (a.e.)%
\begin{equation}\label{eq4.20}
u(t)=\frac{2\sqrt{\frac{\gamma{c_1}-\beta}{\gamma}}}{\ch\left(\pm\sqrt{\beta-\gamma{c_1}}t+a\right)},
\end{equation}

\noindent where
$a=\Arch\left(\frac{2\sqrt{\frac{\beta-\gamma{c_1}}{\gamma}}}{u_B(0)}\right)$.
If $\gamma>0$, then from (\ref{eq4.19}), one can obtain the
candidate optimal control function as follows.%
\begin{itemize}
  \item If $\beta-\gamma{c_1}<0$, then (a.e.)%
\begin{equation}\label{eq4.21}
u(t)=2\sqrt{\frac{\gamma{c_1}-\beta}{\gamma}}\sec\left(\pm\sqrt{\gamma{c_1}-\beta}t+a\right),
\end{equation}

\noindent where
$a=\arcsec\left(\frac{1}{2}\sqrt{\frac{\gamma}{\gamma{c_1}-\beta}}u_B(0)\right)$.

\item If $\beta-\gamma{c_1}>0$, then (a.e.)%
\begin{equation}\label{eq4.22}
u(t)=\frac{2\sqrt{\frac{\beta-\gamma{c_1}}{\gamma}}}{\sh\left(\pm\sqrt{\beta-\gamma{c_1}}t+a\right)},
\end{equation}

\noindent where
$a=\Arsh\left(\frac{2\sqrt{\frac{\beta-\gamma{c_1}}{\gamma}}}{u_B(0)}\right)$.
\end{itemize}

{\it{Remark:}} In the case of $\alpha=0$ and $\gamma\neq0$, the
optimal control function can be expressed by the Weierstrass
elliptic function only when the discriminant $g_2^3-27g_3^2$ of
(\ref{eq4.11}) is nonzero. When $g_2^3-27g_3^2=0$, the optimal
control function is reduced to elementary functions. In comparison
with the case that $g_2^3-27g_3^2\neq0$, the case of
$g_2^3-27g_3^2=0$ occurs with only a probability of zero. The
similar arguments also can be made in the following two cases.

\noindent {\it{2.~Case $\alpha\neq0$ and $\gamma=0$.}}

In this case, control term $B$ of system~(\ref{eq1.1}) is
parabolic. Accordingly, Eq.(\ref{eq4.9}) can be simplified as (a.e.)%
\begin{equation}\label{eq4.25}
(\dot{u})^2=\alpha{u^3}+\beta{u^2}-2\alpha{c_1}u-2\beta{c_1}-2c_2,
\end{equation}

\noindent With $x=\frac{1}{4}\alpha{u}+\frac{1}{12}\beta$,
(\ref{eq4.25}) can be recast as (a.e.)%
\begin{equation}\label{eq4.26}
(\dot{x})^2=4x^3-g_2x-g_3,
\end{equation}

\noindent where $g_2=\frac{1}{12}(\beta^2+6\alpha^2c_1)$ and
$g_3=\frac{1}{216}(18\alpha^2\beta{c_1}+27\alpha^2c_2-\beta^3)$.

If the discriminant of (\ref{eq4.26}) $g_2^3-27g_3^2\neq0$, one
can again obtain the candidate optimal control function in terms
of Weierstrass elliptic function (a.e.)
\begin{equation}\label{eq4.27}
u(t)=\frac{1}{\alpha}\left[4\mathscr{G}(t+a;g_2,g_3)-\frac{\beta}{3}\right],
\end{equation}

\noindent where
$a=\mathscr{G}^{-1}\left(\frac{1}{4}[\alpha{u_B(0)}+\frac{\beta}{3}];g_2,g_3\right)$.

Otherwise, $4x^3-g_2x-g_3$ has repeated zeros. Accordingly, the
candidate optimal control function can be obtained from
(\ref{eq4.25}) or (\ref{eq4.26}) as follows.%
\begin{itemize}
\item If $g_3<0$, then
\begin{equation}\label{eq4.28}
     u(t)=-\frac{6}{\alpha}\sqrt[3]{g_3}\left(\frac{1+\exp(\pm\sqrt{-6\sqrt[3]{g_3}}~t+a)}{1-\exp(\pm\sqrt{-6\sqrt[3]{g_3}}~t+a)}\right)^2
            +\frac{4}{\alpha}\sqrt[3]{g_3}-\frac{\beta}{3\alpha},
\end{equation}

\noindent where
$a=\ln\frac{\sqrt{\frac{1}{4}(\alpha{u_B(0)}+\frac{\beta}{3})-\sqrt[3]{g_3}}-\sqrt{-\frac{3}{2}\sqrt[3]{g_3}}}
            {\sqrt{\frac{1}{4}(\alpha{u_B(0)}+\frac{\beta}{3})-\sqrt[3]{g_3}}+\sqrt{-\frac{3}{2}\sqrt[3]{g_3}}}$.

\item If $g_3=0$, then
\begin{equation}\label{eq4.29}
u(t)=\frac{4}{\alpha(t{\pm}a)^2}-\frac{\beta}{3\alpha},
\end{equation}

\noindent where
$a=\pm\frac{2}{\sqrt{\alpha{u_B(0)}+\frac{\beta}{3}}}$.

\item If $g_3>0$, then
\begin{equation}\label{eq4.30}
u(t)=\frac{6}{\alpha}\sqrt[3]{g_3}\tan\left(\pm\sqrt{\frac{3}{2}\sqrt[3]{g_3}}~t+a\right)
            +\frac{4}{\alpha}\sqrt[3]{g_3}-\frac{\beta}{3\alpha},
\end{equation}

\noindent where
$a=\arctan\sqrt{\frac{\alpha{u_B(0)}+\frac{\beta}{3}-4\sqrt[3]{g_3}}{6\sqrt[3]{g_3}}}$.
\end{itemize}

{\it{3. Case $\alpha\neq0$ and $\gamma\neq0$.}}

In this case, none of the drift term $A$ and the control term $B$
is parabolic. Let%
\begin{equation}\label{eq4.31}
\left\{
\begin{array}{l}
f(x)=\frac{\gamma}{4}x^4+{\alpha}x^3+(\beta-{\gamma}c_1)x^2-2{\alpha}c_1x+\gamma{c_1}^2-2\beta{c_1}-2c_2,\\
g_2=\frac{\gamma}{4}c+\frac{\alpha^2}{2}c_1+\frac{1}{12}(\beta-{\gamma}c_1)^2,\\
g_3=\frac{\gamma(\beta-{\gamma}c_1)}{24}c-\frac{\alpha^2(\beta-{\gamma}c_1)c_1}{24}-\frac{(\beta-{\gamma}c_1)^3}{216}-\frac{\alpha^2{\gamma}c_1^2}{16}-\frac{\alpha^2}{16}c,\\
\end{array}\right.
\end{equation}

\noindent where $c=\gamma{c_1}^2-2\beta{c_1}-2c_2$. Let $x_1$,
$x_2$, $x_3$ and $x_4$ denote the roots of the equation $f(x)=0$.

When $g_2^3-27g_3^2\neq0$, it can be verified that
$x_i\neq{x_j}~(1\leq{i<j}\leq4)$. It is known that the solution of
Eq.(\ref{eq4.10}) still can be written down explicitly in terms of
the Weierstrass function~(see \cite{Whittaker1} chapter XX), which
is given by (a.e.)%
\begin{equation}\label{eq4.34}
u(t)=x_0+\frac{6f'(x_0)}{24\mathscr{G}(t+a;g_2,g_3)-f''(x_0)},
\end{equation}

\noindent where
$a=\mathscr{G}^{-1}\left(\frac{f'(x_0)}{4(u_B(0)-x_0)}+\frac{1}{24}f''(x_0);g_2,g_3\right)$
and $x_0\in\{x_1, x_2, x_3, x_4\}$.

When $g_2^3-27g_3^2=0$, the polynomial $f(x)$ in (\ref{eq4.31})
has repeated zeros. Then, there are four different situations
accordingly.

{\it{1)~$f(x)$ has only one 2-fold zero.}}

In this case, without loss of generality, we can assume that
$x_4~{\neq}~x_1=x_2~{\neq}~x_3~{\neq}~x_4$. Accordingly,
(\ref{eq4.10}) can be rewritten as (a.e.)%
\begin{equation}\label{eq4.35}
(\dot{u})^2=\frac{\gamma}{4}(u-x_1)^2[u^2-(x_3+x_4)u+x_3x_4].
\end{equation}

\noindent From (\ref{eq4.35}), one can compute the
corresponding candidate optimal control function as (a.e.)%
\begin{equation}\label{eq4.36}
u(t)=\frac{2\gamma(x_1-x_3)(x_1-x_4)}{|\gamma(x_3-x_4)|\sin\left(\pm\frac{1}{2}\sqrt{-\gamma(x_1-x_3)(x_1-x_4)}~t+a\right)-\gamma(2x_1-x_3-x_4)}+x_1,
\end{equation}

\noindent where
$a=\arcsin\frac{\gamma(2x_1-x_3-x_4)(u_B(0)-x_1)+2\gamma(x_1-x_3)(x_1-x_4)}
     {(u_B(0)-x_1)|\gamma(x_3-x_4)|}$, when
$\gamma(x_1-x_3)(x_1-x_4)<0$; and%
\begin{equation}\label{eq4.37}
\noindent\begin{array}{l}
u(t)=\frac{\gamma(2x_1-x_3-x_4)+2\left(\exp(\pm\frac{1}{2}\sqrt{\gamma(x_1-x_3)(x_1-x_4)}~t+a)-\frac{\gamma(2x_1-x_3-x_4)}{2\sqrt{\gamma(x_1-x_3)(x_1-x_4)}}\right)\sqrt{\gamma(x_1-x_3)(x_1-x_4)}}
    {\left(\exp(\pm\frac{1}{2}\sqrt{\gamma(x_1-x_3)(x_1-x_4)}~t+a)-\frac{\gamma(2x_1-x_3-x_4)}{2\sqrt{\gamma(x_1-x_3)(x_1-x_4)}}\right)^2-\gamma}+x_1,
\end{array}
\end{equation}

\noindent where $a=\ln\left(
\frac{\sqrt{\gamma{(u_B(0)-x_1)^2}+\gamma(2x_1-x_3-x_4)(u_B(0)-x_1)+\gamma(x_1-x_3)(x_1-x_4)}+\sqrt{\gamma(x_1-x_3)(x_1-x_4)}}{u_B(0)-x_1}\right.$\\

~~~~~~~~~~~~~~~~~~~~~~~~~~~~~~~~~~~~~~~~~~~~~$\left.+\frac{\gamma(2x_1-x_3-x_4)}{2\sqrt{\gamma(x_1-x_3)(x_1-x_4)}}\right)$,

\noindent when $\gamma(x_1-x_3)(x_1-x_4)>0$.

{\it{2)~$f(x)$ has two different 2-fold zeros.}}

In this case, without loss of generality, it can be assumed that
$x_1=x_2~{\neq}~x_3=x_4$. Accordingly,
(\ref{eq4.10}) can be rewritten as (a.e.)%
\begin{equation}\label{eq4.37}
(\dot{u})^2=\frac{\gamma}{4}(u-x_1)^2(u-x_3)^2.
\end{equation}

\noindent If $\gamma>0$ (there is no nontrivial solution for
$\gamma<0$), from (\ref{eq4.37}), we have%
\begin{equation}\label{eq4.38}
u(t)=-\frac{(x_1-x_3)\exp\left(\pm\frac{\gamma}{2}(x_1-x_3)t+a\right)}{\exp\left(\pm\frac{\gamma}{2}(x_1-x_3)t+a\right)-1},
\end{equation}

\noindent where $a=\ln\frac{u_B(0)-x_1}{u_B(0)-x_3}$.

{\it{3)~$f(x)$ has one 3-fold zero.}}

It can be assumed, accordingly, that $x_1=x_2=x_3~{\neq}~x_4$.
Then, (\ref{eq4.10}) can be rewritten as (a.e.)%
\begin{equation}\label{eq4.39}
(\dot{u})^2=\frac{\gamma}{4}(u-x_1)^3(u-x_4).
\end{equation}

\noindent The corresponding candidate optimal control function is
given by (a.e.)%
\begin{equation}\label{eq4.40}
u(t)=\frac{4\gamma(x_1-x_4)}{\gamma^2(x_1-x_4)^2(a\pm\frac{1}{2}t)^2-4\gamma}+x_1,
\end{equation}

\noindent where
$a=-\frac{2}{\gamma(x_1-x_4)(u_b(0)-x_1)}\sqrt{\gamma(u_b(0)-x_1)^2+\gamma(x_1-x_4)(u_b(0)-x_1)}$.

{\it{4)~$f(x)$ has one 4-fold zero.}}

In this case, we can assume that $x_1=x_2=x_3=x_4$. Then,
(\ref{eq4.10}) can be rewritten as (a.e.)%
\begin{equation}\label{eq4.41}
(\dot{u})^2=\frac{\gamma}{4}(u-x_1)^4.
\end{equation}

\noindent If $\gamma>0$ (there is no real solution for
$\gamma<0$), from (\ref{eq4.41}) one can obtain the candidate
optimal control function as (a.e.)%
\begin{equation}\label{eq4.42}
u(t)=\frac{1}{a\pm\frac{\sqrt{\gamma}}{2}t}+x_1,
\end{equation}

\noindent where $a=\frac{1}{u_B(0)-x_1}$.

\section{Examples}\label{sec5}

In this section, we give two examples for illustration.

{\it{Example~1:}} Consider the case when $A=K_x+2K_z$, $B=K_x$
and the target evolution matrix has the form
$X_f=e^{\theta{K_z}}$~($\theta\in\mathbb{R}$).

It can be checked that the abnormal optimal control
$u_\text{abnormal}(t)=-\frac{\left<A,B^\dag\right>}{\left<B,B^\dag\right>}
=-1$ steers system (\ref{eq1.1}) from the initial state $I_2$ to
the final state $X_f$ in time $T_f=(\theta~\text{mod}~4\pi)/2$,
with the performance measure given by%
\begin{equation}\label{eq5.10}
J(u_\text{abnormal})=\int_0^{(\theta~\text{mod}~4\pi)/2}u_\text{abnormal}(t)dt=(\theta~\text{mod}~4\pi)/2.
\end{equation}

\noindent Actually, taking the matrix $S$ as%
\begin{equation}\label{eq5.11}
S=s_zK_z,
\end{equation}

\noindent we have%
\begin{equation}\label{eq5.12}
\begin{array}{l}
   H(S;\lambda_o;u;X_o(t))=\left<S,{X^o}(t)^{-1}[A+u(t)B]{X^o}(t)\right>\\
   ~~~~~~~~~~=\left<s_zK_z,e^{-2tK_z}[K_x+2K_z+uK_x]e^{2tK_z}\right>\\
   ~~~~~~~~~~=\left<s_zK_z,(1+u)[\cos(2t)K_x-\sin(2t)K_y]+2K_z\right>\\
   ~~~~~~~~~~\equiv2s_z,\\
\end{array}
\end{equation}

\noindent which is minimized with respect to $u$ by the abnormal
optimal control $u_\text{abnormal}(t)$.\\

{\it{Example~2:}} Suppose that the system (\ref{eq1.1}) is given
by $A=K_z$ and $B=-K_x+K_y$. Consider the optimal steering problem
with the terminal condition $X_f=e^{-2K_x+2K_y}$.

Write the matrix $S$ as%
\begin{equation}\label{eq5.1}
S=s_xK_x+s_yK_y+s_zK_s.
\end{equation}

\noindent In order to achieve target evolution matrix optimally,
according to the Eqs. (\ref{eq4.3a}) and (\ref{eq4.3b}), the
coefficients $s_x$, $s_y$ and $s_z$ in (\ref{eq5.1}) are required
to satisfy%
\begin{equation}\label{eq5.2}
s_z=\frac{\sinh(2\sqrt{2})}{\sqrt{2}[\cosh(2\sqrt{2})-1]}(s_x+s_y).
\end{equation}

\noindent Fig~\ref{fig1} shows the resulting evolutions.

\noindent\begin{center}%
\noindent\begin{figure}[h]%
\centering%
\includegraphics[totalheight=2.7in]{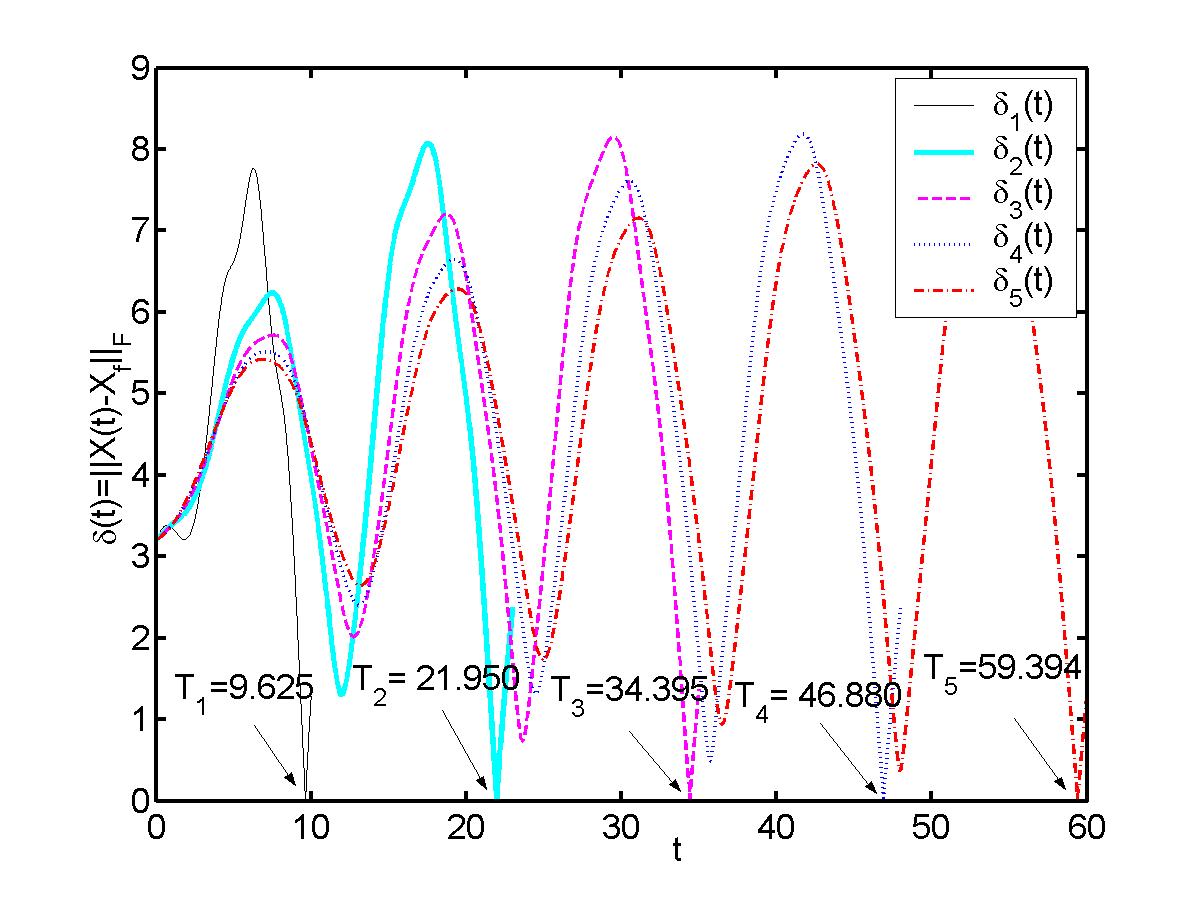}%
\vspace{-0.6cm}%
\caption{The Distance Between the State $X(t)$ and the Target
$X_f$ at time $t$. A numerical treatment shows that when the two
independent parameters $(s_x,s_y)$ take the values of
$(-2.09895,0.99801)$, $(-0.844738,0.406815)$,
$(-0.523766,0.254595)$, $(-0.379568,0.185591)$,
$(-0.297326,0.145922)$, $\cdots$, the target evolution matrix
$X_f=e^{-2K_x+2K_y}$ can be reached at time $T_1=9.625$,
$T_2=21.950$, $T_3=34.395$, $T_4=46.880$, $T_5=59.394$, $\cdots$.
Where the norm $||X||_F$ is defined as
$||X||_F=\sqrt{\sum_{i,j}|X_{ij}|^2}$.}\label{fig1}
\end{figure}
\end{center}
\noindent Correspondingly, the control functions given by%
\begin{equation}\label{eq5.3}%
\begin{array}{l}
  |u_1(t)|=\sqrt{2\mathscr{G}(t+0.5740-1.2502i;9.8362,4.4871)+2.4990} \\
  ~~~~~~~~~~=\frac{1.0547}{\sqrt{\mathscr{G}(t+0.5740;9.8362,4.4871)+1.2495}},
\end{array}
\end{equation}
\begin{equation}\label{eq5.4}%
\begin{array}{l}
   |u_2(t)|=\sqrt{2\mathscr{G}(t+0.6705-2.0412i;3.6531,1.3108)+1.2418}\\
   ~~~~~~~~~~=\frac{0.3487}{\sqrt{\mathscr{G}(t+0.6705;3.6531,1.3108)+0.6209}},
\end{array}
\end{equation}
\begin{equation}\label{eq5.5}%
\begin{array}{l}
   |u_3(t)|=\sqrt{2\mathscr{G}(t+0.7051-2.5221i;2.6390,0.8203)+0.9956}\\
   ~~~~~~~~~~=\frac{0.2045}{\sqrt{\mathscr{G}(t+0.7051;2.6390,0.8203)+0.4978}},
\end{array}
\end{equation}
\begin{equation}\label{eq5.6}%
\begin{array}{l}
   |u_4(t)|=\sqrt{2\mathscr{G}(t+0.7247-2.8643i;2.2391,0.6435)+0.8954}\\
   ~~~~~~~~~~=\frac{0.1442}{\sqrt{\mathscr{G}(t+0.7247;2.2391,0.6435)+0.4477}},
\end{array}
\end{equation}
\begin{equation}\label{eq5.7}%
\begin{array}{l}
   |u_5(t)|=\sqrt{2\mathscr{G}(t+0.7370-3.1118i;2.0254,0.5543)+0.8416}\\
   ~~~~~~~~~~=\frac{0.1115}{\sqrt{\mathscr{G}(t+0.7370;2.0254,0.5543)+0.4208}},
\end{array}
\end{equation}
\begin{equation}\nonumber%
   \vdots
\end{equation}
\noindent as shown in Fig~\ref{fig2}, are candidate optima.
\noindent\begin{center}%
\noindent\begin{figure}[h]%
\centering%
\includegraphics[totalheight=2.7in]{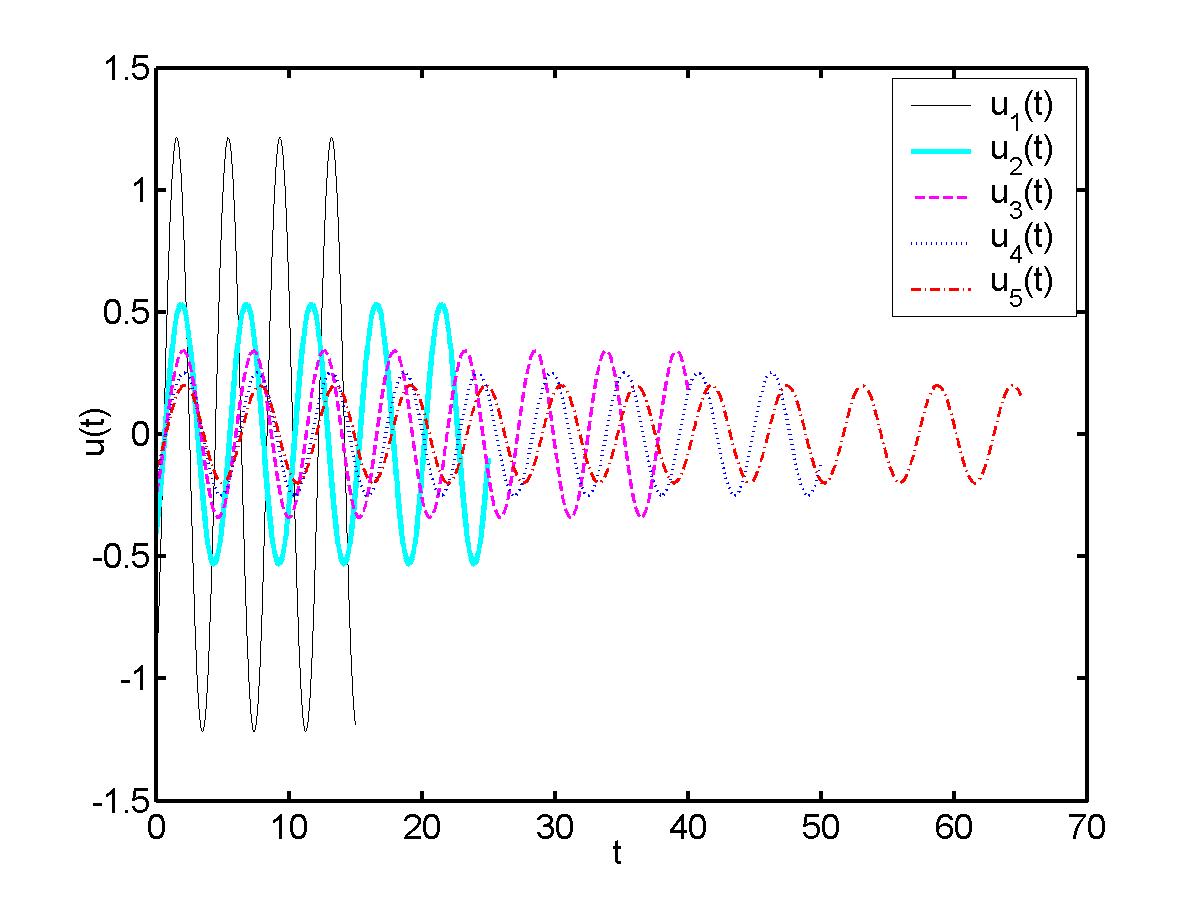}%
\vspace{-0.6cm}%
\caption{The Candidate Optimal Control Functions.}\label{fig2}
\end{figure}
\end{center}\noindent
The corresponding performance measures are shown in
Fig~\ref{fig3}. It is easy to observe that there is a tradeoff
between the consumed time $T$ and the cost index $J(u)$.
\noindent\begin{center}%
\noindent\begin{figure}[h]%
\centering%
\includegraphics[totalheight=2.7in]{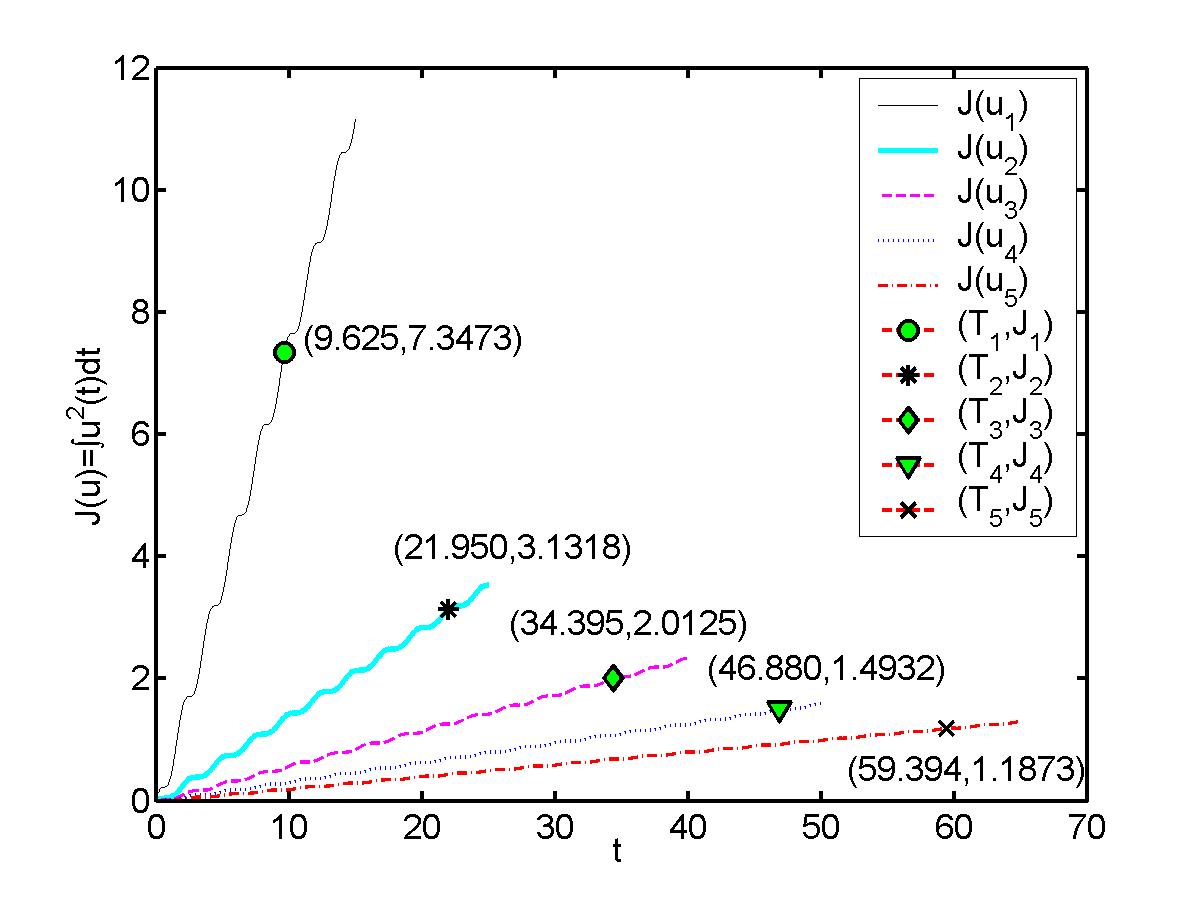}%
\vspace{-0.6cm}%
\caption{The Performance Measures Corresponding to the Candidate
Optimal Controls. Where $J(u_1)=7.3473$, $J(u_2)=3.1318$,
$J(u_3)=2.0125$, $J(u_4)=1.4932$, $J(u_5)=1.1873$,
$\cdots$.}\label{fig3}
\end{figure}
\end{center}\noindent

According to the decomposition algorithm given in~\cite{wjw1},
piecewise constant control laws can be designed to achieve the
target evolution $X_f=e^{-2K_x+2K_y}$ as well. One possible
control law is given by%
\begin{equation}\label{eq5.8}%
\bar{u}(t)=
\begin{array}{l}
\end{array}\left\{%
\begin{array}{ll}
    0, & \hbox{when $t\in\left[0,t_1+2n_1\pi\right)$;} \\
    \frac{c}{\sqrt{2}}, & \hbox{when $t\in\left[t_1+2n_1\pi,t_1+t_2+2n_1\pi\right]$;} \\
    0, & \hbox{when $t\in\left(t_1+t_2+2n_1\pi,2t_1+t_2+2(n_1+n_2)\pi\right]$,} \\
\end{array}%
\right.
\end{equation}

\noindent where $n_1,n_2\in{\mathbb{N}^+}$, $c>1$ and%
\begin{equation}\nonumber%
\left\{\begin{array}{l}
   t_1=2\arccot\left(-c\cth\sqrt{2}\right.-\left.\sqrt{c^2\cth^2\sqrt{2}-1}\right),\\
   t_2=\frac{2}{\sqrt{c^2-1}}\Arcth\left(\frac{1}{\sqrt{c^2-1}}\sqrt{c^2\cth^2\sqrt{2}-1}\right).\\
\end{array}\right.
\end{equation}

\noindent The corresponding performance measure is given by%
\begin{equation}\label{eq5.9}%
J(\bar{u})=\frac{c^2}{\sqrt{c^2-1}}\Arcth\left(\frac{1}{\sqrt{c^2-1}}\sqrt{c^2\cth^2\sqrt{2}-1}\right).
\end{equation}

\noindent It can be verified that $J(\bar{u})$ increases
monotonously with the increase of $c$. Since%
\begin{equation}\label{eq5.9}%
\lim\limits_{c\rightarrow1}J(\bar{u})=\frac{1}{\sqrt{\cth^2\sqrt{2}-1}}\approx1.9351,
\end{equation}

\noindent as a comparison, the performance measures $J(u_4)$ and
$J(u_5)$ are approximately $13$ and $39$ percent, respectively,
less than $J(\bar{u})$.

\section{Conclusion}\label{sec6}

In order to minimize the decoherence effect, energy optimal
control problem for the quantum systems evolving on the noncompact
Lie group $SU(1,1)$ is taken into account in this paper. We showed
that explicit expressions for the optimal control functions can be
obtained analytically. To minimize the considered quadratic
performance measure, the control functions with respect to
abnormal extremals are constant functions of time $t$, while those
with respect to normal extremals are expressed by the Weierstrass
elliptic function.

\section*{Acknowledgment}

This research was supported in part by the National Natural
Science Foundation of China under Grant Nos 60674039, 60433050 and
60635040. Tzyh-Jong~Tarn would also like to acknowledge partial
support from the U.S. Army Research Office under Grant
W911NF-04-1-0386.

The authors would like to thank Dr. Re-Bing Wu for his helpful
suggestions.

\end{document}